\documentclass[fleqn,twoside,twocolumn,nofootinbib,showkeys]{revtex4} 
\usepackage[sec,nocpr]{ujp} 

\begin{document}
\title[Generalization of Polynomial Invariants]
{GENERALIZATION OF POLYNOMIAL\\ INVARIANTS AND HOLOGRAPHIC PRINCIPLE\\ FOR KNOTS AND LINKS}%
\author{A.M. Pavlyuk}
\affiliation{\bitp}
\address{\bitpaddr}
\email{pavlyuk@bitp.kiev.ua}
\udk{539.12; 517.984}  \razd{\secix}

\autorcol{A.M.\hspace*{0.7mm}Pavlyuk}

\newcommand{\ve}{\varepsilon}
\newcommand{\wt}{\widetilde}
\newcommand{\cl}{\centerline}
\newcommand{\ts}{\textstyle}

\setcounter{page}{673}%

\begin{abstract}
We formulate the holographic principle for knots and links. For
the ``space'' of all knots and links, torus knots $T({\it
2m+1\!,2})$ and torus links $L({\it 2m,2})$
 play the role
of the ``boundary'' of this space. Using the holographic principle,
we
 find the skein relation of knots and links with the help of
the recurrence relation for polynomial invariants of torus knots
$T({\it 2m+1\!,2})$ and torus links $L({\it 2m,2})$. As an example
of the application of this principle, we derive the Jones skein
relation and its generalization with the help of some variants of
$(q,p)$-numbers, related with $(q,p)$-deformed bosonic
oscillators.
\end{abstract}
 \maketitle

\section{Introduction}
The interaction of knot theory and physics manifests itself from
the knotting of physics to the physicalness of knots~\cite{Ka1,Ka2,At}.
The construction of invariants of knots and links on the basis of deep
physical ideas shows the physical nature of knots and links.
In particular, the applications of the Jones polynomials and its generalizations~\cite{Wi}
demonstrate that knots may play an important role
in the unification of general relativity and quantum physics~\cite{Ba1,Ba2}.
Witten introduced a heuristic definition of the Jones polynomial invariants
 in terms of topological quantum field theory.

A powerful modern tool for physical investigations, the holographic principle (HP),
 was introduced in 1993 by 't Hooft~\cite{Ho},
 when studying the quantum mechanical features of black holes at the Planck scale.
It claims that the total information contained in a volume of space
is encoded on the boundary of this space. From the 't Hooft conjecture,
it follows that our world at Planckian energies can be described as
2d, rather than 3d; e.g., 3d quantum gravity can be described
by 2d topological field theory. In 1997, the HP was used by
Maldacena~\cite{Ma} to develop AdS/CFT correspondence, which
manifests a holographic duality between guantum gravity in the 5d
anti-de Sitter space and the conformal field theory on its 4d boundary.
Verlinde~\cite{Ve} recently proposed that the HP provides a natural
mechanism for the Newton gravity and the Einstein gravity to emerge due to
entropic forces.

Thus, the HP allows one to deal with the theory in a D-dimensional
space instead of studying a (D+1)-dimensional theory. Among
possible generalizations of the HP, one can seek such ones that
describe the properties of the whole system on the basis of the
properties of its part, which can be considered as a ``boundary''.
In this paper, we make a conjecture on the existence of the HP for
the ``space'' of all knots and links. As the prompting point of
this conjecture, we can use the fact that the transformations of
3d knot and link structures are described with the help of the
transformations of their 2d planar diagrams. To express the HP for
knots and links mathematically, we consider the torus knots
$T(2m+1,2)$ and the torus links $L(2m,2)$ with $m$ to be a
non-negative integer, as playing the role of the ``boundary'' for
the ``space'' of all knots and links.
We will show that such an important information (concerning all
knots and links), as the skein relation, can be obtained from the
``boundary'' information -- the recurrence relation for the
polynomial invariants of torus knots $T(2m+1,2)$ and torus links
$L(2m,2)$, which, in turn, can be obtained  with the help of the
properties of $(q,p)$-numbers~\cite{GP1,Pa2}. The ``boundary''
position of these torus knots and links together with the fact
that a lot of information concerning all knots and links can be
``decoded''\ from the ``boundary information'', allows us to speak
about the HP for knots and links. Using this HP, we will find the
Jones skein relation and its 2-variable generalization.

\section{Holographic Principle\\ for Knots and Links}

Polynomial invariants belong to the most important characteristics
of knots and links. They are defined with the help of the skein
relation, which looks as
\begin{equation}\label{skein}
P_{L_{+}}(t)=l_{1}(t)P_{L_{O}}(t)+l_{2}(t)P_{L_{-}}(t),
\end{equation}
where $l_{1}(t),\, l_{2}(t)$ are the coefficients. The skein
relation~(\ref{skein}) connects three polynomials $P_{L_{+}}(t),\,
P_{L_{O}}(t),$ and  $P_{L_{-}}(t),$ corresponding to underscript
knots (links). Here, $L_{+}$ denotes an arbitrary knot or link
(all knots and links are also called as ``Links''). From the
initial Link $L_{+}$ by the so-called surgery operation (namely,
the elimination of a crossing),  we obtain the Link $L_{O}$. From
the same initial Link $L_{+}$ by another surgery operation
switching off a crossing, we obtain another Link $L_{-}$.
For~(\ref{skein}) to be the skein relation, it must satisfy the
Reidemeister moves. The skein relation~(\ref{skein}) together with
the normalization condition for the unknot,
\begin{equation}
\label{unknot} P_{\rm unknot}=1\,,
\end{equation}
put a definite polynomial in correspondence to every knot and
link.

Therefore, to operate with the skein relation~(\ref{skein}), it is
necessary to know the coefficients $\,l_{1}\,,\, l_{2}\,.$ The
Alexander skein relation follows from~(\ref{skein}), if
\[
l_{1}=t^{1/2}-t^{-{1/2}},\quad l_{2}=1.
\]
For the Jones skein relation,
\[
 l_{1}=t(t^{1/2}-t^{-{1/2}}),\quad  l_{2}=t^{2}.
\]
The 2-variable HOMFLY skein relation corresponds to
\[ l_{1}=tz,\quad l_{2}=t^{2}.\]

Let us consider both torus knots and links denoted as $L_{n,2}$,
$n=0,1,2, \mbox{...},$ including torus knots $T(2m+1,2)$ and torus
links $L(2m,2),\ m=0,1,2, \mbox{...}\,.$ The surgery operation of
elimination turns torus knot/link  $L_{n, 2}$  into  torus
link/knot $L_{n-1, 2}$, and the surgery operation of switching
turns torus knot/link $L_{n, 2}$ into torus knot/link $L_{n-2,
2}.$

Therefore, the following
recurrence relation for torus knots and links $L_{n, 2}$ can be obtained from~(\ref{skein}):
\[ P_{L_{n+1, 2}}(t)=l_{1}(t)P_{L_{n, 2}}(t)+l_{2}(t)P_{L_{n-1, 2}}(t)\,.\]
Rewriting it in a simpler form
\begin{equation}
\label{skein1} P_{n+1, 2}(t)=l_{1}P_{n, 2}(t)+l_{2}P_{n-1, 2}(t)\,
\end{equation}
and comparing with~(\ref{skein}), we see that there exists the
correspondence
\begin{equation}
\begin{array}{l}\label{cor} P_{L_{+}}(t) {\leftrightarrow}  P_{n+1,
2}(t),\quad P_{L_{O}}(t){\leftrightarrow} P_{n, 2}(t),
\\[3mm]
P_{L_{-}}(t){\leftrightarrow} P_{n-1, 2}(t),
\end{array}
\end{equation}
which allows one, having~(\ref{skein1}), to write the skein
relation~(\ref{skein}). For this, it is enough to take $l_{1}$ and
$l_{2}$ from~(\ref{skein1}) and to put them into~(\ref{skein}).
For the verification, we use the Reidemeister moves. This is an
example of the HP for knots and links.
For the ``space'' of all knots and links, the torus knots and
links $L_{n,2}$ are considered as the ``boundary'' of this space.

From the recurrence relation~(\ref{skein1}) for the class of torus
knots and links $L_{n,2}$, one finds the recurrence relation for
its subclass. For example, only for torus knots $T(2m+1,2)$ (or
only for torus links $L(2m,2)$) (which allows one to find the
skein relation~(\ref{skein1}) as well)~\cite{GP1}, we have
 \begin{equation}\label{skein2}
P_{n+2, 2}(t)=k_{1}P_{n, 2}(t)+k_{2}P_{n-2, 2}(t)\,,
\end{equation}
where
\begin{equation}\label{k} k_{1}=l_{1}^{2}+2l_{2},\quad k_{2}=-l_{2}^{2}.
\end{equation}
With the help of the normalization condition~(\ref{unknot}), we
also find
\begin{equation}
\label{32} P_{1,2}=1,\quad P_{3,2}=k_{1}+k_{2}.
\end{equation}

\section{Jones Skein Relation\\ from the Holographic Principle}

The $(q,p)$-number (or the so-called structural function for a
$(q,p)$-deformed bosonic oscillator) is defined as~\cite{CJ}
\begin{equation}
\label{qp-def} [n]_{q,p}={{\textstyle {q^{n}-p^{n}}}\over{\textstyle
{q-p}}}.
\end{equation}
The recurrence relation for $(q,p)$-numbers~(\ref{qp-def})\linebreak
looks as
\begin{equation} \label{qp-rec}
[n+1]_{q,p}=(q+p)[n]_{q,p}-qp[n-1]_{q,p}.
\end{equation}

Let us make the following substitution in~(\ref{qp-def}):
\begin{equation}
\label{qpt} q\rightarrow q^{3},\quad p\rightarrow q.
\end{equation}
It reduces the $(q,p)$-numbers to the $(q^{3},q)$-numbers
\begin{equation}
\label{q3q} [n]_{q^{3},q}={{\textstyle
{q^{3n}-q^{n}}}\over{\textstyle {q^{3}-q}}},
\end{equation}
which satisfy the recurrence relation
\begin{equation}\label{q3q-rec}
[n+1]_{q^{3},q}=(q+q^{3})[n]_{q^{3},q}-q^{4}[n-1]_{q^{3},q}.
\end{equation}
With the help of the three-step algorithm~\cite{GP1,Pa2}
and~(\ref{q3q-rec}), we can find the Jones skein relation.

First, we introduce the polynomials $J_{n,2}(q)$ satisfying the
recurrence relation repeating~(\ref{q3q-rec}) (i.e., the
recurrence relation only for torus knots $T(2m+1,2)$ or only for
torus links $L(2m,2)$):
\begin{equation}
\label{jones-skein2}
J_{n+2,2}(q)=(q+q^{3})J_{n,2}(q)-q^{4}J_{n-2,2}(q).
\end{equation}
 From~(\ref{32}), we obtain
\begin{equation}
\label{J13} J_{1,2}(q)=1,\quad J_{3,2}(q)=q+q^{3}-q^{4}.
\end{equation}

Second,  we formulate the recurrence relation for all polynomials
$J_{n, 2}(q)$, which leads to the corresponding skein relation.
From (\ref{skein2}) and (\ref{jones-skein2}), we have $k_{1}=
q+q^{3}$, $k_{2}=-q^{4}$. Using (\ref{k}), one finds
\begin{equation}
\label{Al}
 l_{1}=q(q^{1/2}-q^{-{1/2}}), \quad l_{2}=q^{2}.
\end{equation}
Therefore, we have the recurrence relation for torus knots
$T(2m+1,2)$ and for torus links $L(2m,2)$, taken together:
\begin{equation}
\label{jones-skein1} J_{n+1, 2}(q)=q(q^{1/2}-q^{-{1/2}})J_{n,
2}(q)+q^{2}J_{n-1, 2}(q)\,.
\end{equation}
From~(\ref{jones-skein1}) in correspondence with the HP for knots
and links in the form~(\ref{cor}), we obtain the Jones skein
relation~\cite{Jo}
\begin{equation}
\label{jones-skein}
J_{+}(q)=q(q^{1/2}-q^{-{1/2}})J_{O}(q)+q^{2}J_{-}(q).
\end{equation}

Third, we find the expression for the Jones polynomials of torus
knots~$T(2m+1,2)$ in terms of the $(q^{3},q)$-numbers:
\begin{equation}
\label{jones-gavr}
J_{2m+1,2}(q)=b_{1}(q)[m+1]_{q^{3},q}-b_{2}(q)[m]_{q^{3},q}.
\end{equation}
Substituting ~(\ref{J13}) into~(\ref{jones-gavr}), we find
$b_{1}=1$ and $b_{2}=q^{4}.$ Thus,
\begin{equation}\label{jones-gavr2}
J_{2m+1,2}(q)=[m+1]_{q^{3},q}-q^{4}[m]_{q^{3},q}.
\end{equation}
It is easy to verify that formula~(\ref{jones-gavr2}) can be
derived from the known general formula for the Jones polynomials
of torus knots,
\begin{equation}\label{jones-torus} J_{n,k}(q)=q^{{(n-1)(k-1)}\over{2}}
{{1-q^{n+1}-q^{k+1}+q^{n+k}}\over {1-q^{2}}},
\end{equation}
if $k=2$.

\section{Generalized Jones Polynomials}

Let us make the substitution
\begin{equation}
\label{q3p} q\rightarrow q^{3}
\end{equation}
in (\ref{qp-def}). From $(q^{3},p)$-numbers
\begin{equation}
\label{q3p-def} [n]_{q^{3},p}={{\textstyle
{q^{3n}-p^{n}}}\over{\textstyle {q^{3}-p}}},
\end{equation}
using the three-step procedure as in the previous section, we
obtain the 2-parameter generalized Jones skein relation:
\begin{equation}
\label{gJ}
 J_{+}(q,p)=(q^{3/2}-p^{1/2})J_{O}(q,p)+q^{3/2}p^{1/2}J_{-}(q,p).
\end{equation}\vspace*{-5mm}

\noindent Formula (\ref{gJ}) can be rewritten as
\[
q^{-{3/4}}p^{-{1/4}} J_{+}(q,p)-q^{3/4}p^{1/4}J_{-}(q,p)=
\]\vspace*{-5mm}
\begin{equation}
\label{gJ-}= (q^{3/4}p^{-{1/4}}-q^{-{3/4}}p^{1/4})J_{O}(q,p).
\end{equation}
By putting $p=q,$ this generalized Jones skein relation
turns into the usual Jones skein relation. 
At last, we have the analytical formula for the simplest torus
knots:
\begin{equation}
\label{AA2m-} J_{2m+1,2}(q,p)=[m+1]_{q^{3},p}-q^{3}p[m]_{q^{3},p}.
\end{equation}
It turns into formula~(\ref{jones-gavr2}) for the Jones polynomial
invariants if $p=q.$

Comparing (\ref{gJ-}) with the skein relation for the HOMFLY
polynomial invariants~\cite{FY}
\begin{equation}
\label{homfly-skein} a^{-1}H_{+}(a, z)-aH_{-}(a, z)=zH_{O}(a, z),
\end{equation} we find that the substitution
\begin{equation}\label{compar}
 a=q^{3/4}p^{1/4}, \quad z=q^{3/4}p^{-{1/4}}-q^{-{3/4}}p^{1/4}
 \end{equation}
turns the 2-parameter generalized Jones skein relation into the
HOMFLY one.\vspace*{1mm}

\section{Concluding Remarks}

In this paper, we state that there exists the holographic
principle for knots and links. We have presented a simple example
that demonstrates obtaining
 the skein relation of knots and links~(\ref{skein}) from
the recurrence relation~(\ref{skein1}) (or~(\ref{skein2})) for the
simplest torus knots $T(2m+1,2)$ and torus links $L(2m,2)$ due to
the direct correspondence between their coefficients  $l_{1}$ and
$l_{2}.$ The formulas expressing this correspondence depend on the
used surgery operations in the definition of \mbox{skein
relation.}

We note that the derivation of the formulas for the polynomial
invariants of arbitrary knots and links on the basis of the known
formulas for torus knots $T(2m+1,2)$ and torus links $L(2m,2)$ would
support the validity of the HP for knots and links. As the first
step in this direction, it is necessary to solve the simpler problem
concerning arbitrary torus knots and links. In fact, we need the key
for decoding the ``boundary information'' and for turning it into
the ``volume information''.

As one more example of the HP for knots and links, we mention the
expansion of an arbitrary polynomial invariant into the sum of
polynomial invariants of torus knots $T(2m+1,2)$ and torus links
$L(2m,2).$ An attempt to do this for torus knots $T(n,3)$ and
$T(n,4)$ was made in~\cite{GP2,Pa1}.


\vspace*{-5mm}
\rezume{%
А.М. Павлюк } {УЗАГАЛЬНЕННЯ ПОЛІНОМІАЛЬНИХ\\ ІНВАРІАНТІВ І
ГОЛОГРАФІЧНИЙ ПРИНЦИП\\ ДЛЯ ВУЗЛІВ І ЗАЧЕПЛЕНЬ} {Сформульовано
голографічний принцип для вузлів і зачеплень. Для ``простору'' всіх
вузлів і зачеплень торичні вузли $T(2m+1,2)$ і торичні зачеплення
 $L(2m,2)$ відіграють роль ``межі'' цього простору.
 Використовуючи голографічний принцип для вузлів і зачеплень,
 знаходимо skein-співвідношення для вузлів
 і зачеплень за допомогою рекурентних співвідношень для поліноміальних інваріантів
 торичних вузлів  $T(2m+1,2)$ і торичних зачеплень
 $L(2m,2).$
Як приклад застосування цього принципу, отримано skein-співвідношення
Джонса і його узагальнення за допомогою певних варіантів
$(q,p)$-чисел, що пов'язані з $(q,p)$-деформованими бозонними
осциляторами.
 }%



\end{document}